\documentclass[12pt,letterpaper]{article}
\usepackage[T1]{fontenc}
\usepackage[utf8]{inputenc}
\usepackage[english]{babel}
\usepackage[margin=0.75in]{geometry}
\usepackage{lmodern}
\usepackage{mathtools,amsfonts,amssymb}
\usepackage[numbers,super,sort&compress]{natbib}
\usepackage[linesnumbered,ruled]{algorithm2e}
\usepackage{hyperref}
\usepackage{authblk}
\usepackage{comment}
\allowdisplaybreaks

\title{Wavefront reconstruction of discontinuous phase objects from optical deflectometry}

\author[1]{Ricardo Legarda-Saenz\footnote{Corresponding author: \texttt{rlegarda@correo.uady.mx} URL: \url{http://clir-lab.org}}}
\author[2]{Jorge L. Flores}
\author[3]{Manuel Servin}
\author[1]{Anabel Martin-Gonzalez}
\affil[1]{CLIR at Facultad de Matematicas, Universidad Autonoma de Yucatan. 97205 Merida, Mexico.}
\affil[2]{Departamento de Electronica, Universidad de Guadalajara. 44840 Guadalajara, Mexico.}
\affil[3]{Centro de Investigaciones en Optica, A. C. 37150 Leon, Mexico.}
\date{\today}

\begin{document}
\maketitle

\begin{abstract}
One of the challenges in phase measuring deflectometry is to retrieve the wavefront from objects that present discontinuities or non-differentiable gradient fields. Here, we propose the integration of such gradients fields based on an $L^{p}$-norm minimization problem. The solution of this problem results in a nonlinear partial differential equation, which can be solved with a fast and well-known numerical methods and doesn’t depend on external parameters. Numerical reconstructions on both synthetic and experimental data are presented that demonstrate the capability of the proposed method.
\end{abstract}

\section{Introduction}
\begin{figure}[ht]
  \begin{center}
   \includegraphics[width=\columnwidth]{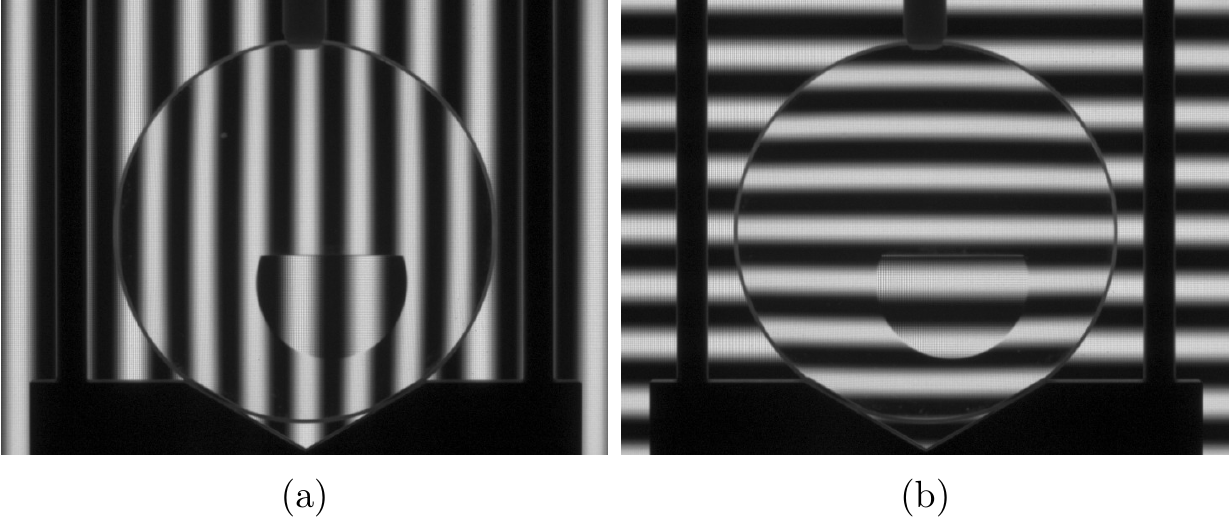}
    \caption{The fringe pattern captured by the camera with (a) $\mathbf{v}_{1} = (1,0)$, and (b) $\mathbf{v}_{2} = (0,1)$.}\label{fig:Franjas}
  \end{center}
\end{figure}

Phase-Measuring Deflectometry is a well established technique for measuring specular surfaces or transparent objects, because of its qualities of non-contact, low-cost, and high-resolution.~\cite{Massig1999,Kammel2008,Huang2018,Xu2020a,Wang2021a,Gao2022a,Gao2022}  For the particular case of measuring transparent phase objects, the typical experimental setup consists of a LCD monitor and CCD camera, ~\cite{Massig1999,Canabal2002,Legarda-Saenz2007a,Vargas2010,Flores2015a} which was originally proposed in reference \citenum{Massig1999}: Assume that the LCD monitor displays a sinusoidal intensity fringe patterns parallel to the y-direction and an ideal transparent phase object is placed in front of it, at a given distance $D$ from it. The acquired fringe pattern will be distorted by the phase object due of it changes the optical path lengths traveled by light rays through it, e.g., if the phase is inhomogeneous in the x-direction, the rays will be deflected by $\alpha D = \partial\phi_{\mathbf{x}}/\partial x$ where $\mathbf{x} = (x,y)$ and the term $\phi_{\mathbf{x}}$ is the optical path length accumulated by a ray traveling through the phase object at the position $\mathbf{x}.$ In other words, the fringe pattern distortion is related with the gradient field of the phase object. Sketches of this experimental setup can be seen in references \citenum{Massig1999} and \citenum{Flores2015a}. An example of captured fringe patterns using this experimental setup are shown in Figure \ref{fig:Franjas}. 

The fringe pattern captured by the camera is described by
\begin{equation}\label{eq:franjas}
I_{\mathbf{x}} = a_{\mathbf{x}} + b_{\mathbf{x}}\cos\left( \frac{2\pi}{q}\left[\mathbf{x}\cdot\mathbf{v} + D\nabla\phi_{\mathbf{x}}\cdot\mathbf{v}\right] \right),  
\end{equation}
where $a_{\mathbf{x}}$ is the background illumination, $b_{\mathbf{x}}$ is the amplitude modulation, $\mathbf{v} = (\cos\varphi,\sin\varphi)$ is the normal direction vector of the pattern displayed on the screen, the term $\nabla\phi_{\mathbf{x}} = (\partial\phi_{\mathbf{x}}/\partial x,\partial\phi_{\mathbf{x}}/\partial y) $ is their wavefront gradient of $\phi_{\mathbf{x}}$, and $q$ is the pattern period displayed on the screen. The resultant measurement of this technique is the directional derivative of the wavefront $\psi$ oriented in the direction of $\mathbf{v}$. 

A common procedure to estimate the wavefront is to acquire two orthogonal directional derivatives and integrate them.~\cite{Fried1977,Hudgin1977a,Frankot1988,Roddier1991,Zou2000} Typically, one choose two orthogonal directions, like $\mathbf{v}_{1} = (1,0)$ and $\mathbf{v}_{2} = (2,0).$ For that, two sequence of fringe patterns, in horizontal and vertical directions, are acquired, where the demodulated phase term is given by 
\begin{align*}
\Psi^{x} &= D\nabla\phi_{\mathbf{x}}\cdot\mathbf{v}_1 \approx \partial\phi_{\mathbf{x}}/\partial x,\\
\Psi^{y} &= D\nabla\phi_{\mathbf{x}}\cdot\mathbf{v}_2 \approx \partial\phi_{\mathbf{x}}/\partial y, 
\end{align*}
from which it will be integrated to obtain the gradient field of the wavefront.~\cite{Fried1977,Hudgin1977a,Frankot1988,Roddier1991}

One of the major actual challenges is to retrieve with high accuracy the wavefront introduced by objects that present discontinuities in their surface, for instance a bifocal lens. In literature there are some methods which care about discontinuities:~\cite{Karacali2003,Reddy2009,DiMartino2018,Queau2018a,Queau2018} Queau et al. present an extensive review of integration methods~\cite{Queau2018a} and propose some variational approaches for integration surfaces with discontinuities.~\cite{Queau2018} As it is pointed out, their variational methods are suitable to recover discontinuities, but they have slow iterative solvers, and they need to properly tune some parameters to successfully recover discontinuities.

In this paper, we present a different approach to integrate the obtained information from a deflectometry setup based on $L^{p}$-norm minimization problem.  The solution of this problem result on a nonlinear partial differential equation (PDE), which can be solved with a fast and well-known numerical methods, with the advantage that this procedure does not depend on external parameters. The outline of this paper is as follows: In Section \ref{sec:metodo} we review the $L^{p}$-norm minimization problem to be solved and present the first-order optimality conditions consisting on a nonlinear PDE. In Section \ref{sec:numerical} we review the numerical solution of the PDE and propose an algorithm to solve it, which results in a simple one. The numerical results on both synthetic and experimental data are presented in Section \ref{sec:experiment} and our conclusions are given in Section \ref{sec:conclusion}.

\section{Methodology}\label{sec:metodo}
Following the solution given in references \citenum{Ghiglia1996} and \citenum{Ghiglia1998} for phase unwrapping, we propose to integrate the wavefront $\phi$ given the measured gradient field $\Psi = \left( \Psi^{x}, \Psi^{y}\right),$ requiring that the derivatives of the wavefront $\phi$ are agreed with the measured gradient field $\Psi$ in the minimum $L^{p}$-norm sense. That is, the integrated wavefront $\phi$ that minimizes
\begin{equation}\label{eq:funcional}
\underset{\phi}{\min}\; E(\phi_{\mathbf{x}}) = \iint_{\Omega}\left\lvert\frac{\partial\phi_{\mathbf{x}}}{\partial x} - \Psi_{\mathbf{x}}^{x}\right\rvert^{p} 
+ \left\lvert\frac{\partial\phi_{\mathbf{x}}}{\partial y} -\Psi_{\mathbf{x}}^{y}\right\rvert^{p}\;d\mathbf{x}
\end{equation}
is the minimum $L^{p}$-norm solution, and $\Omega \subseteq \mathbb{R}^{2}$ is the domain of integration.

To obtain the solution of the problem expressed in Equation  (\ref{eq:funcional}), the first-order optimality condition or Euler-Lagrange equation has to be derived, resulting in the following nonlinear PDE
\begin{equation}\label{eq:gradiente}
- \;\frac{\partial}{\partial x}\left[\left(\frac{\partial\phi_{\mathbf{x}}}{\partial x} - \Psi_{\mathbf{x}}^{x}\right)
\left\lvert\frac{\partial\phi_{\mathbf{x}}}{\partial x} - \Psi_{\mathbf{x}}^{x}\right\rvert^{p-2} \right] 
- \;\frac{\partial}{\partial y}\left[\left(\frac{\partial\phi_{\mathbf{x}}}{\partial y} - \Psi_{\mathbf{x}}^{y}\right)
\left\lvert\frac{\partial\phi_{\mathbf{x}}}{\partial y} - \Psi_{\mathbf{x}}^{y}\right\rvert^{p-2}\right] = 0,
\end{equation}
with boundary conditions
\[\;\left[ \left(\frac{\partial\phi_{\mathbf{x}}}{\partial x} - \Psi_{\mathbf{x}}^{x}\right)
\left\lvert\frac{\partial\phi_{\mathbf{x}}}{\partial x} - \Psi_{\mathbf{x}}^{x}\right\rvert^{p-2}, \; 
\left(\frac{\partial\phi_{\mathbf{x}}}{\partial y} - \Psi_{\mathbf{x}}^{y}\right)
\left\lvert\frac{\partial\phi_{\mathbf{x}}}{\partial y} - \Psi_{\mathbf{x}}^{y}\right\rvert^{p-2}\right] \cdot \hat{\mathbf{n}} = 0,\]
where $\hat{\mathbf{n}}$ denotes the unit outer normal vector to the boundary. The constant $p$ is removed to avoid cancellation when one sets $p = 0.$

Using the following substitution 
\begin{equation}\label{eq:UV}
U_{\mathbf{x}} = \left\lvert\frac{\partial\phi_{\mathbf{x}}}{\partial x} - \Psi_{\mathbf{x}}^{x}\right\rvert^{p-2},\qquad 
V_{\mathbf{x}} = \left\lvert\frac{\partial\phi_{\mathbf{x}}}{\partial y} - \Psi_{\mathbf{x}}^{y}\right\rvert^{p-2},
\end{equation}
Equation (\ref{eq:gradiente}) becomes 
\begin{equation}\label{eq:gradienteUV}
- \;\frac{\partial}{\partial x}\left[\left(\frac{\partial\phi_{\mathbf{x}}}{\partial x} - \Psi_{\mathbf{x}}^{x}\right) U_{\mathbf{x}} \right] 
- \;\frac{\partial}{\partial y}\left[\left(\frac{\partial\phi_{\mathbf{x}}}{\partial y} - \Psi_{\mathbf{x}}^{y}\right) V_{\mathbf{x}} \right] = 0,
\end{equation}
with boundary conditions
\[\;\left[ \left(\frac{\partial\phi_{\mathbf{x}}}{\partial x} - \Psi_{\mathbf{x}}^{x}\right) U_{\mathbf{x}}, \; 
\left(\frac{\partial\phi_{\mathbf{x}}}{\partial y} - \Psi_{\mathbf{x}}^{y}\right) V_{\mathbf{x}} \right] \cdot \hat{\mathbf{n}} = 0.\]

In principle, the functional proposed in Eq. (\ref{eq:gradienteUV}) can take any value of $p$. However, there are values that allow us to integrate surfaces with discontinuities or non differentiable gradient fields. For values $p \geq  2,$ the solution tends to be smooth and the discontinuities of the gradient field will be lost, which is also undesirable for our approach. On the contrary, for values $p < 2,$ particularly $p \leq 1,$ the solution tends to follow the gradient field of the input data;~\cite{Bloomfield1983,Scales1988} that is, the minimum $L^{p}$-norm finds a solution whose gradients agree best with those of the input data in the sense of $L^{p}$-norm. As it was reported in references \citenum{Ghiglia1996} and \citenum{Ghiglia1998}, this approach provides well-behaved computations and the solution deviates from the gradient field of the input data in a very few places. One important point of the solution shown in Eq. (\ref{eq:gradienteUV}) is that generates its own data-dependent weights when $p \neq 2.$ These weights allow, on the one hand, to do the adjustment of least squares in continuous regions where the residuals are relatively small, associated with noise; on the other hand, they do not allow over-penalization in regions where the residuals are large, which are associated with regions with discontinuities or non differentiable gradient fields.

\section{Numerical solution}\label{sec:numerical}
Let $u_{i,j} = u(x_i,y_j)$ to denote the value of a function $u_{\mathbf{x}}$ at point $(x_i,y_j)$ defined on $\Omega = [a,b]\times[c,d],$ where the sampling points are $x_i = a + (i-1)h_x,\quad y_j = c + (j-1)h_y,$ with $1 \leq i \leq M,\; 1 \leq j \leq N,$ $h_x = (b-a)/(M-1),\; h_y = (d-c)/(N-1)$ and $M, N$ are the number of points in the discrete grid. We use $u$ to represent any of the variables $\phi,\Psi^{x},\Psi^{y}, U, V$ defined in the previous equations. Derivatives are approximated using standard forward and backward finite difference schemes \[\delta_{x}^{\pm}u_{i,j} = \pm\frac{u_{i,j\pm 1}-u_{i,j}}{h_x}\quad\text{and}\quad\delta_{y}^{\pm}u_{i,j} = \pm\frac{u_{i\pm 1j}-u_{i,j}}{h_y}.\]  The gradient and the divergence are approximated as \[\triangledown  u_{i,j} = (\delta_{x}^{+}u_{i,j},\delta_{y}^{+}u_{i,j})\quad\text{and}\quad\triangledown\cdot\triangledown  u_{i,j} = \delta_{x}^{-} (\delta_{x}^{+}u_{i,j}) + \delta_{y}^{-}(\delta_{y}^{+}u_{i,j}),\] respectively. 

Hence the numerical approximation of the Equation (\ref{eq:gradienteUV}) is given by 
\begin{equation}\label{eq:gradienteUVdisc}
\begin{aligned}
&- \;\delta_{x}^{-}\left[\left(\delta_{x}^{+}\phi_{i,j} - \Psi_{i,j}^{x}\right) U_{\mathbf{x}} \right] 
- \;\delta_{y}^{-}\left[\left(\delta_{y}^{+}\phi_{i,j} - \Psi_{i,j}^{y}\right) V_{\mathbf{x}} \right] = \\
&\left(\phi_{i,j} - \phi_{i,j-1} - \Psi_{i,j-1}^{x}\right)U_{i,j-1} - \left(\phi_{i,j+1} - \phi_{i,j} - \Psi_{i,j}^{x}\right)U_{i,j}\\
  + &\left(\phi_{i,j} - \phi_{i-1,j} - \Psi_{i-1,j}^{y}\right)V_{i-1,j} - \left(\phi_{i+1,j} - \phi_{i,j}  - \Psi_{i,j}^{y}\right)V_{i,j}
  = 0,
\end{aligned}
\end{equation}
with boundary conditions
\[\;\left[ \left(\phi_{i,j+1} - \phi_{i,j} - \Psi_{i,j}^{x}\right) U_{i,j}, \;\left(\phi_{i+1,j} - \phi_{i,j}  - \Psi_{i,j}^{y}\right) V_{i,j}\right] \cdot \hat{\mathbf{n}} = 0,\]  where for simplicity and without loss of generality we consider $h_x = h_y  = 1.$

The stability and convergence of the numerical solution are better behave when the values of the terms $U_{i,j}$ and $V_{i,j}$ lie in range (0,1),~\cite{Bloomfield1983,Scales1988,Ghiglia1996} so we rewrite Equation (\ref{eq:UV}) as 
\begin{equation}\label{eq:UVnorm}
\begin{aligned}
U_{i,j} &= \left\lvert\phi_{i,j+1}-\phi_{i,j} - \Psi_{i,j}^{x}\right\rvert^{p-2}
= \frac{\epsilon }{\left\lvert\phi_{i,j+1}-\phi_{i,j} - \Psi_{i,j}^{x}\right\rvert^{2-p} + \epsilon},\\
V_{i,j} &= \left\lvert\phi_{i+1,j}-\phi_{i,j} - \Psi_{i,j}^{y}\right\rvert^{p-2}
= \frac{\epsilon }{\left\lvert\phi_{i+1,j}-\phi_{i,j} - \Psi_{i,j}^{y}\right\rvert^{2-p} + \epsilon},
\end{aligned}
\end{equation}
where $\epsilon$ is a value to force the normalization of the terms $U_{i,j}$ and $V_{i,j}$ into the range (0, 1). In our experiments, we found that $\epsilon = 0.1$ fulfill this requirement.

Equation (\ref{eq:gradienteUV}) and their discrete version Eq. (\ref{eq:gradienteUVdisc})  are a nonlinear PDE because of the terms U and V are functions of the wavefront $\phi$ and the measured gradient field $\Psi.$ Here, an iterative method will be adopted to solve Eq. (\ref{eq:gradienteUVdisc})~\cite{Bloomfield1983,Scales1988}: first, given an initial value $\phi,$ the terms $U$ and $V$ are computed; then, the terms $U$ and $V$ are held fixed and Eq. (\ref{eq:gradienteUVdisc}) is solved using preconditioned conjugate gradient~\cite{Shewchuk1994,Golub1996}. With the current solution $\phi,$ the terms $U$ and $V$ are updated and a new solution is computed. This process is repeated until convergence. 

Now, we present the algorithm to solve the discretization of Eq. (\ref{eq:gradienteUV}). First we arrange Eq. (\ref{eq:gradienteUVdisc}) in  matrix form $\mathbf{A}\phi = \mathbf{b}$ as
\begin{equation}\label{eq:Axb}
  \begin{aligned}
  &-\left(\phi_{i,j+1}U_{i,j} + \phi_{i+1,j}V_{i,j} + \phi_{i,j-1}U_{i,j-1} + \phi_{i-1,j}V_{i-1,j}\right)\\
  &+\left(U_{i,j-1} + U_{i,j} + V_{i-1,j} + V_{i,j}\right)\phi_{i,j}\\
  &= \Psi_{i,j-1}^{x}U_{i,j-1} - \Psi_{i,j}^{x}U_{i,j} + \Psi_{i-1,j}^{y}V_{i-1,j} - \Psi_{i,j}^{y}V_{i,j}.
  \end{aligned}
\end{equation}
Notice that $\mathbf{A}$ is a $MN\times MN$ sparse matrix which depends on the terms $U$ and $V,$ so it needs to be constructed at each iteration.  The resultant linear system is symmetric, semi-positive definite and weakly diagonally dominant, so we use the the preconditioned conjugate gradient (PCG) to solve this system using the implementation proposed in Ref. \citenum{Shewchuk1994}. The explicit structure  of the algorithm is given in Appendix \ref{sec:apendice}, where we use the incomplete LU factorization with no fill-in as preconditioning technique to estimate matrix $\mathbf{M}.$~\cite{VanderVorst2003,Saad2003,Chen2005a,Legarda-Saenz2022}

\section{Numerical experiments}\label{sec:experiment}
To illustrate the performance of the proposed technique, we carried out two sets of numerical experiments using a Intel\textsuperscript{\tiny\textregistered} Core\textsuperscript{\tiny\texttrademark} i7 @ 2.40 GHz laptop with Debian GNU/Linux 11 (bullseye) 64-bit and 16 GB of memory. For our experiments, we programmed all the functions using C language, GNU g++ 10.2 compiler and Intel\textsuperscript{\tiny\textregistered} oneAPI Math Kernel Library, 2022.1 release.\footnote{\url{https://www.intel.com/content/www/us/en/developer/tools/oneapi/onemkl.html}}  The values used as stopping criteria in Algorithm \ref{algoritmo} were $ k_{max} = 100,$ $tol = 10^{-3},$ $l_{max} = 1.5MN,$ and $\kappa = 0.005.$

\begin{figure}[ht]
  \begin{center}
   \includegraphics[width=\columnwidth]{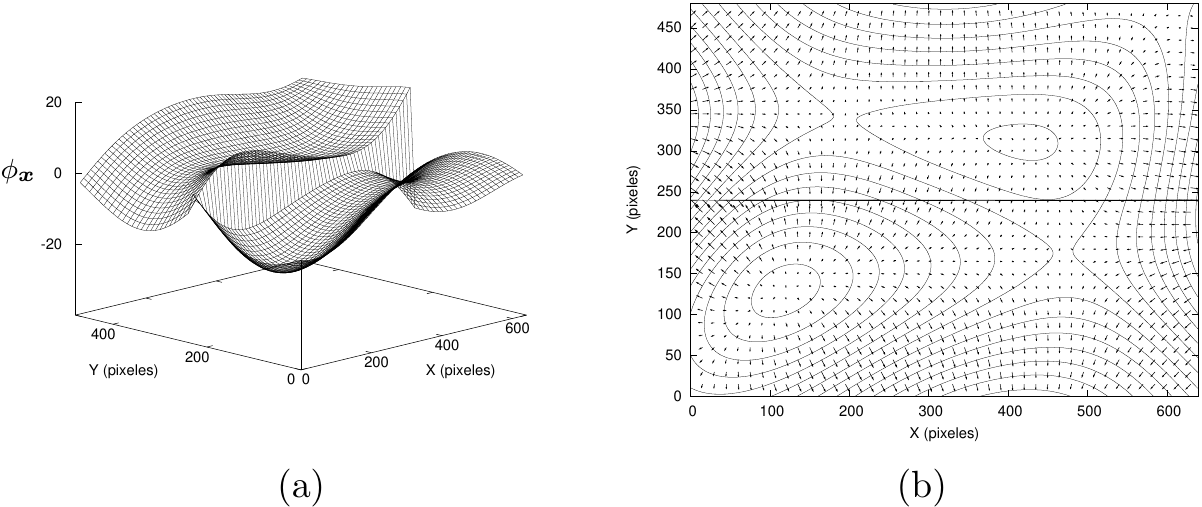}
    \caption{(a) Synthetic wavefront and (b) its gradient field.}\label{fig:SinteticoOriginal}
  \end{center}
\end{figure}

Our first experiment was the integration of the gradient field obtained from a synthetic wavefront defined as
\begin{equation}\label{eq:sintetico}
	\begin{aligned}
		\Theta_{\mathbf{x}}  &= 15(1-x)^{2}\exp\left[-x^{2} - (y+1)^2\right]\ldots\\  &- 50*\left(x/5 - x^{3} - y^{5}\right)\exp\left[-x^{2} - y^2\right]\ldots\\ 
		&- 5/3\exp\left[-(x+1)^{2} - y^2\right],\\
		\phi_{\mathbf{x}} &=
		\begin{cases}
			\Theta_{\mathbf{x}}       & \quad \text{if } x \geq 0\\
			-\Theta_{\mathbf{x}}       & \quad \text{if } x < 0
		\end{cases},
	\end{aligned}
\end{equation}
evaluated in the domain $\Omega = [-1,1]\times[-1,1],$ using M = 480 and N = 640 as the number of points in the discrete grid. Figure \ref{fig:SinteticoOriginal} shows this synthetic wavefront and their gradient field.

We made several estimations of this synthetic wavefront using different values of $p.$ For each values of $p,$ we corrupted the gradient field shown in Figure \ref{fig:SinteticoOriginal}, panel (b) with different levels of Gaussian noise. We use a normalized error $Q$ to compare the resultant estimations; this error is defined as~\cite{Perlin2016} \[Q\left(\mu,\nu\right) = \frac{\| \mu - \nu\|_{2}}{\| \mu\|_{2} + \|\nu\|_{2}},\] where $\mu$ and $\nu$ are the signals to be compared. The normalized error values vary between zero (for perfect agreement) and one (for total disagreement). 

\begin{figure}[ht]
  \begin{center}
   \includegraphics[width=\columnwidth]{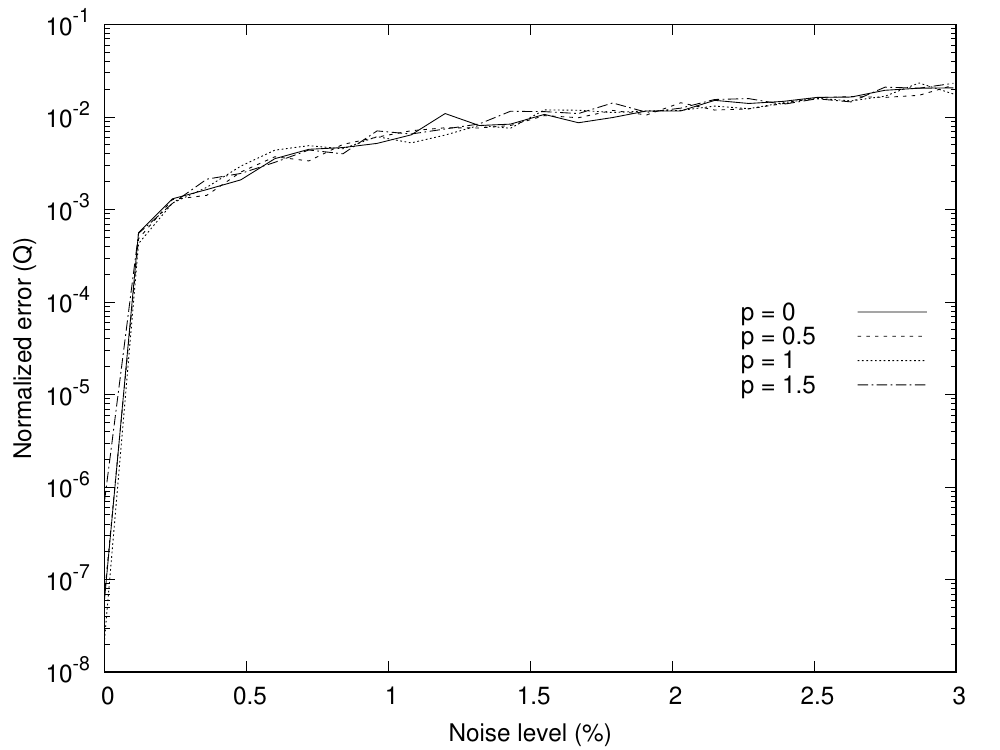}
    \caption{Normalized error for different noise levels.}\label{fig:Error}
  \end{center}
\end{figure}

\begin{figure}[ht]
  \begin{center}
   \includegraphics[width=\columnwidth]{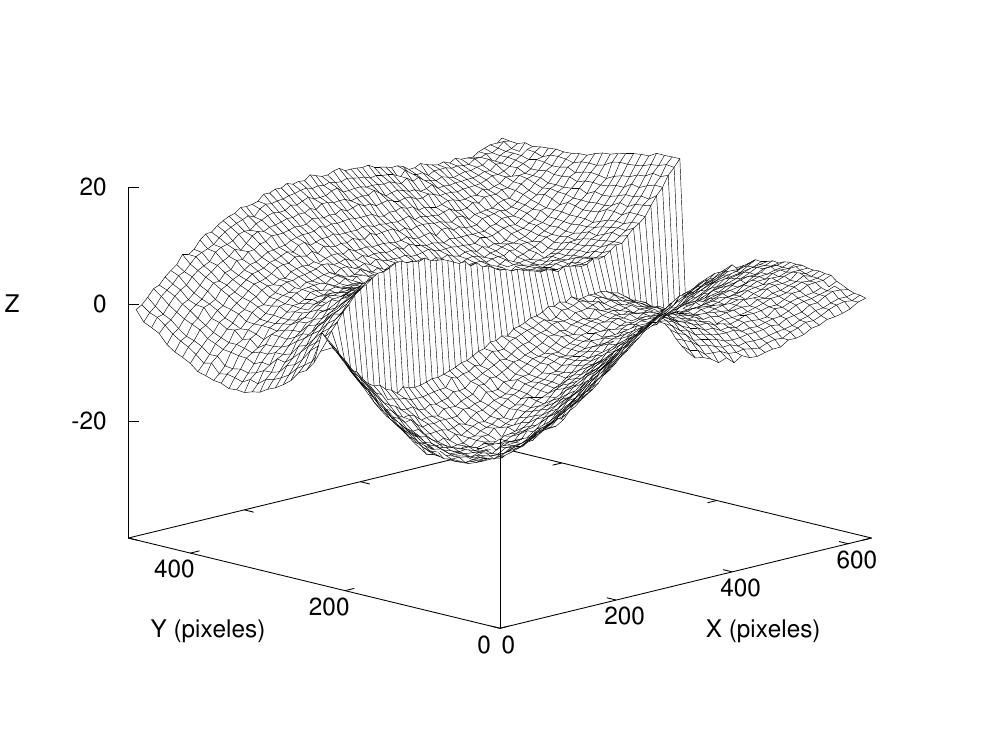}
    \caption{Estimated wavefront from synthetic information using $p$ = 0.}\label{fig:SinteticoEstimado}
  \end{center}
\end{figure}

Figure \ref{fig:Error} shows the errors obtained for the different levels of noise and different values of $p.$ Each integration experiment takes approximately 1400 iterations, about 18 seconds to perform each experiment on the previously mentioned computer equipment. An example of the resultant reconstruction using $p = 0$ is shown in Figure \ref{fig:SinteticoEstimado}, where the gradient field was corrupted with 3\% of Gaussian noise. In addition, Table \ref{tab:tabla} shows the number of iterations and normalized error $Q$ using different values of $p$ in the synthetic wavefront estimation without noise.
\begin{table}[h]
	\centering 
	\caption{Obtained values for synthetic wavefront estimations without noise}\label{tab:tabla}
	\begin{tabular}{ c c c }
		\hline
		p& iterations&Q\\
		\hline
		0 & 1390 & 2.5e-08 \\
		0.5 &1292 & 2.7e-08 \\
		1.0 &1388 & 1.7e-08 \\
		1.5 &1023 & 1.4e-06
	\end{tabular}
\end{table}

\begin{figure}[ht]
	\begin{center}
		\includegraphics[width=\columnwidth]{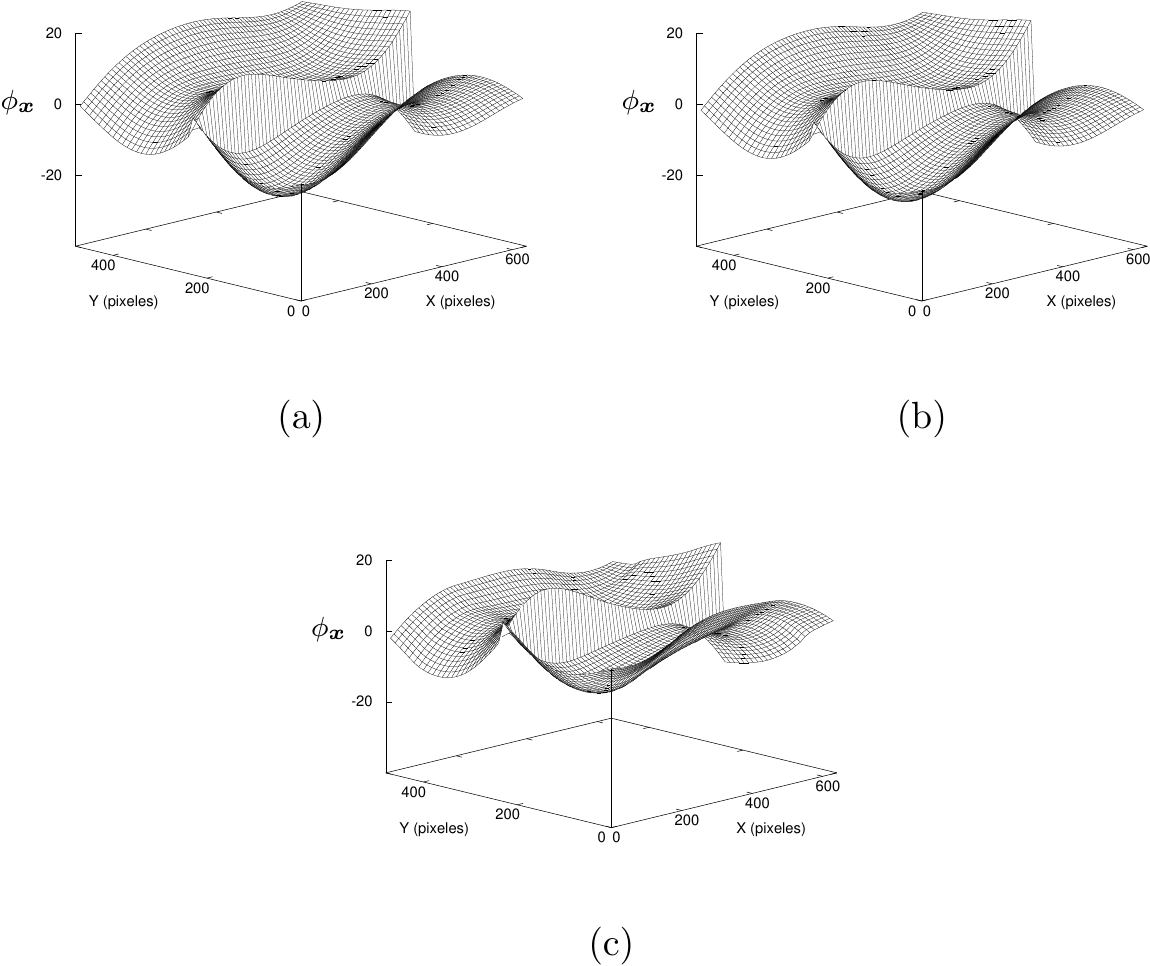}
		\caption{Estimated wavefront from synthetic information using: (a) $p = 1$,  (b) reference \citenum{Queau2018}, and (c) $p = 2$.}\label{fig:SinteticoEstimado1}
	\end{center}
\end{figure}

\begin{figure}[ht]
	\begin{center}
		\includegraphics[width=\columnwidth]{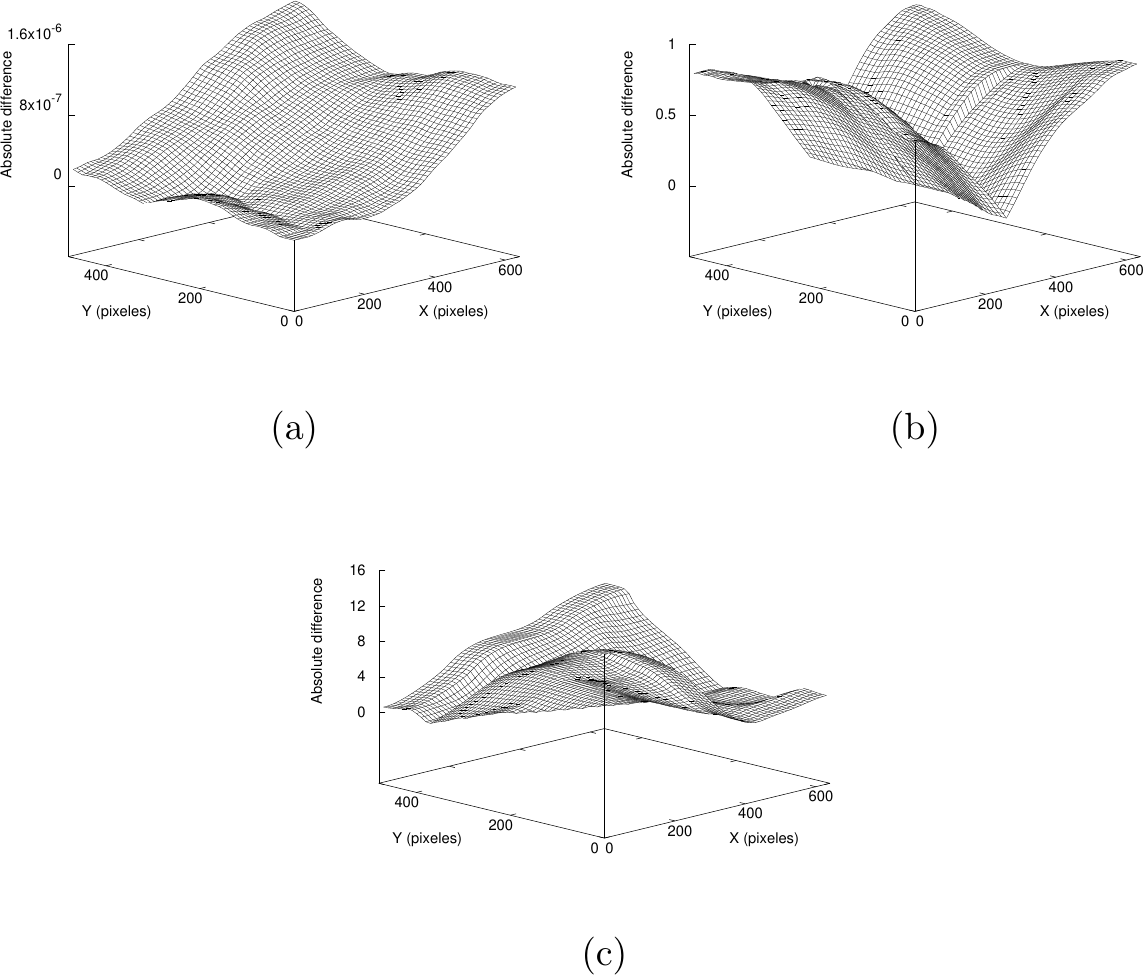}
		\caption{Absolute difference obtained from wavefront estimation of the synthetic information using: (a) $p = 1$,  (b) reference \citenum{Queau2018}, and (c) $p = 2$.}\label{fig:SinteticoEstimado2}
	\end{center}
\end{figure}

As can be observed, the integration based on Eq. (\ref{eq:gradiente}), and their implementation shown in Algorithm \ref{algoritmo}, successfully reconstruct the wavefront $\phi_{\mathbf{x}}.$ As expected, the normalized error grows as noise is added to the gradient field; however, the proposed solution performs the reconstruction without affecting the scale of the wavefront and preserving the edges; all this for different values of $p$, using similar processing times. From Table \ref{tab:tabla} and Figure \ref{fig:Error}, we find that the $p$ values can be defined between $0 \le p \le 1$ without affecting the reconstruction performance. This is consistent with what has been described in classical references such as \citenum{Bloomfield1983}, \citenum{Scales1988} and \citenum{Ghiglia1998}. In order to compare the performance of our proposal, a second experiment was the gradient field integration without noise of the synthetic wavefront defined in Eq. (\ref{eq:sintetico}), using $p = 1$, $p = 2$, and the variational approach described in section 4.2 of reference \citenum{Queau2018}. The estimated wavefronts are shown in Figure \ref{fig:SinteticoEstimado1} and the absolute differences are shown in Figure \ref{fig:SinteticoEstimado2}. As it can be observed, the integration based on Eq. (\ref{eq:gradiente}) shows a better performance in the wavefront reconstruction, preserving the discontinuities.

\begin{figure}[ht]
  \begin{center}
   \includegraphics[width=\columnwidth]{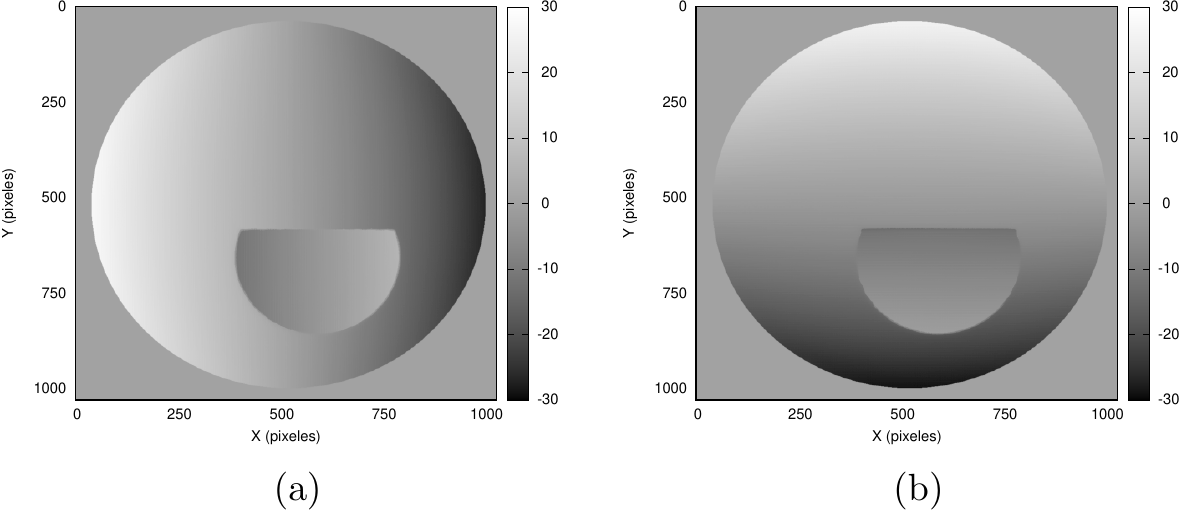}
    \caption{Experimental gradient field: (a) $\Psi^{x}$ and (b) $\Psi^{y}$.}\label{fig:ExperimentalGradientes}
  \end{center}
\end{figure}

\begin{figure}[ht]
  \begin{center}
   \includegraphics[width=\columnwidth]{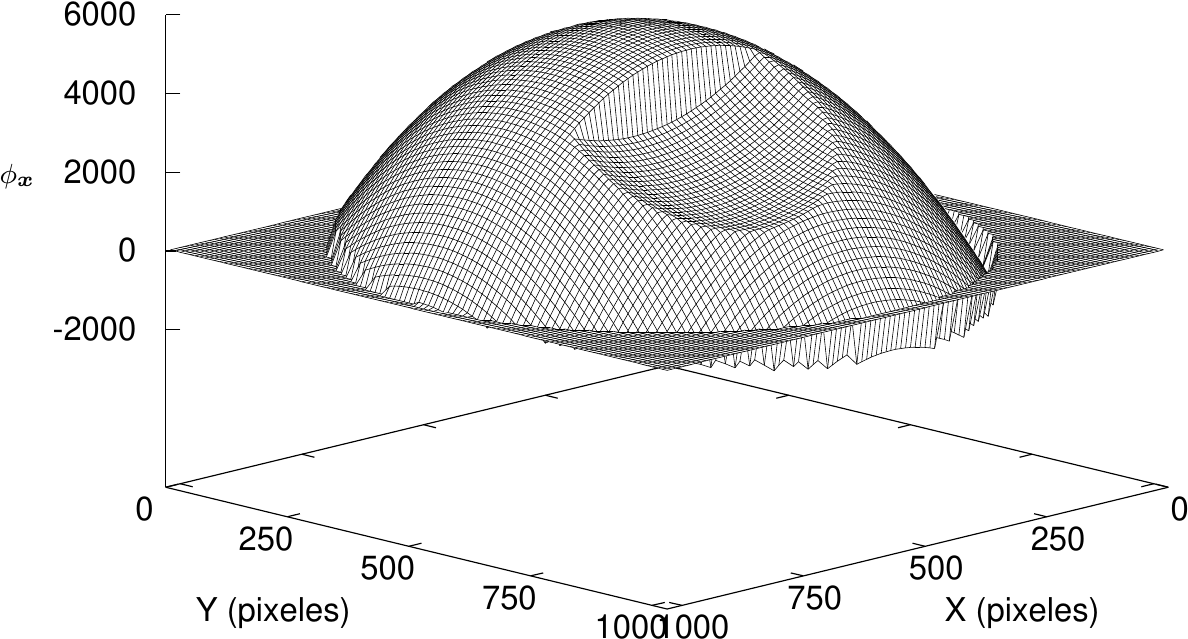}
    \caption{Estimated wavefront from experimental information using $p$ = 0.}\label{fig:ExperimentalEstimado}
  \end{center}
\end{figure}

The third experiment was the integration of the gradient field shown in Figure \ref{fig:ExperimentalGradientes}, measured from a bifocal lens using the technique described in the introduction. Example of the fringe pattern and the experimental setup can be observed in Figure \ref{fig:Franjas}. The estimated wavefront is shown in Figure \ref{fig:ExperimentalEstimado} and the time employed to obtain this estimation was 447 seconds approximately using 8950 iterations with $p = 0$. In this experiment we show the performance of the proposed method on the processing of experimental information. As can be observed from the wavefront reconstruction of the bifocal lens, the proposed solution adequately integrates the two surfaces that are on the lens and its edge with respect to the background.

\section{Discussion of results and conclusion}\label{sec:conclusion}
The presented integration method allows recovering the wavefront of objects that present discontinuities in their surface or non differentiable gradient fields from phase measuring deflectometry. The numerical solution of Eq. (\ref{eq:gradiente}) results in a very simple algorithm, as it can be observed in Algorithm \ref{algoritmo}. An extra advantage of the proposed solution is its feasibility to be implemented with different techniques like multigrid or to be parallelized. This will be one aim of our future research.

\appendix
\section{Appendix}\label{sec:apendice}
In this appendix, we show the proposed algorithm to solve Eq. (\ref{eq:Axb}).
\begin{algorithm}
	\DontPrintSemicolon
	\KwData{the measured gradient field $\Psi = \left( \Psi^{x}, \Psi^{y}\right)$ and $p < 2.$}
	\KwResult{the integrated wavefront $\phi.$}
	\BlankLine
	$k \leftarrow 0,\;error \leftarrow 1$\;
	$\phi_{i,j}^{k}\leftarrow\text{random values}$\;
	\While{$k < k_{max}\;\mathbf{ and}\; error > tol $}
	{
		compute $U$ and $V$, Eq. (\ref{eq:UVnorm})\;
		solve Eq. (\ref{eq:Axb}) using PCG:\;
		\Begin
		{
			construct $\mathbf{A}$ and $\mathbf{b}$ \;
			estimate preconditioning matrix $\mathbf{M}$ from $\mathbf{A}$\;\label{matrizM}
			$\phi^{k+1} \leftarrow\phi^{k}$\;
			$l \leftarrow 0$\;
			$\mathbf{r}\leftarrow\mathbf{b} - \mathbf{A}\phi^{k+1}$\;
			$\mathbf{d}\leftarrow\mathbf{M}^{-1}\mathbf{r}$\;
			$\delta_{new}\leftarrow\mathbf{r}^{T}\mathbf{d}$\;
			$\delta_{0}\leftarrow\delta_{new}$\;
			\While{ $l < l_{max}\;\mathbf{ and}\; \delta_{new} > \kappa^{2}\delta_{0}$}
			{
				$\mathbf{q}\leftarrow\mathbf{A}\mathbf{d}$\;
				$\alpha\leftarrow\delta_{new} / \mathbf{d}^{T}\mathbf{q}$\;
				$\phi^{k+1} \leftarrow\phi^{k+1} + \alpha\mathbf{d}$\;
				\eIf{ $l$ is divisible by $\sqrt{MN}$ }
				{$\mathbf{r}\leftarrow\mathbf{b} - \mathbf{A}\phi^{k+1}$\;}
				{$\mathbf{r}\leftarrow\mathbf{r} - \alpha\mathbf{q}$\;} 
				$\mathbf{s}\leftarrow\mathbf{M}^{-1}\mathbf{r}$\;
				$\delta_{old}\leftarrow\delta_{new}$\;
				$\delta_{new}\leftarrow\mathbf{r}^{T}\mathbf{s}$\;
				$\beta\leftarrow\delta_{new} / \delta_{old}$\;
				$\mathbf{d}\leftarrow\mathbf{s} + \beta\mathbf{d}$\;
				$l\leftarrow l +1$\;  
			}
		}
		$error\leftarrow\|\phi_{i,j}^{k+1}-\phi_{i,j}^{k}\| / \|\phi_{i,j}^{k}\|$\;    
		$k\leftarrow k +1$\;  
	}
	\caption{\(L^p\)-norm integration algorithm.}\label{algoritmo}
\end{algorithm}

\newpage
\bibliographystyle{unsrt}    

\newpage

\end{document}